\documentclass[12pt]{elsarticle}

\usepackage{amssymb,amsmath,amsthm,amscd,amsfonts}

\newcommand{\diffspec}{\operatorname{Spec^\Delta}}
\newcommand{\poly}[2]{#1^{\{#2\}}}
\newcommand{\height}{\operatorname{ht}}

\newcommand{\diffht}{\operatorname{ht^{\Delta}}}
\newcommand{\type}{\operatorname{type^\Delta}}
\newcommand{\diffdim}{\operatorname{dim^\Delta}}

\newtheorem{theorem}{Theorem}
\newtheorem{lemma}[theorem]{Lemma}
\newtheorem{proposition}[theorem]{Proposition}
\newtheorem{corollary}[theorem]{Corollary}

\newtheorem*{statement*}{Statement}
\newtheorem*{theorem*}{Theorem}
\newtheorem*{lemma*}{Lemma}

\theoremstyle{definition}
\newtheorem{definition}[theorem]{Definition}
\newtheorem*{definition*}{Definition}
\newtheorem{example}{Example}
\newtheorem*{example*}{Example}

\theoremstyle{remark}
\newtheorem{remark}[theorem]{Remark}
\newtheorem*{remark*}{Remark}

\begin{document}

\begin{frontmatter}

\title{Differential Krull dimension in differential polynomial extensions}
\author{Ilya Smirnov\fnref{authorlabel}}
\ead{ilya.smrnv@gmail.com}
\fntext[authorlabel]{The work of the author was partially supported by the NSF Grant CCF-1016608.}

\address{Department of Mechanics and Mathematics, Moscow State University, Russia, 119991}

\begin{abstract}
In the present paper, we investigate the differential Krull
dimension of rings of differential polynomials. We prove a
differential analogue of the Special Chain Theorem and show that
some important classes of differential rings have no anomaly of
differential Krull dimension.
\end{abstract}


\end{frontmatter}

\section{Introduction}

In the theory of algebraic differential equations, we have the
notion of the differential dimension of the set of solutions.
Roughly speaking, the dimension shows how many independent solutions
we have. But this characteristic describes the ``global behavior'',
that is the behavior of a ``typical'' solution. If we study
solutions near a given one, many different anomalies occur. The
examples can be found in~\cite{Rit1,Rit2}.

In~\cite[p.~607]{SelKl}, it is noted that the notion of differential
Krull dimension should play the desired role locally. And one of the
most important questions is whether a differential integral domain
being differentially finitely generated over a field is
differentially catenary~\cite[p.~608]{SelKl}. This question is
related with the behavior of descending chains of differential prime
ideals. Such chains were investigated in~\cite{Ros} and~\cite{J1}.
But these papers did not provide a full answer. The last deep result
about differential Krull dimension was obtained in~\cite{J1} and,
after this paper, there were no results in this direction because of
the lack of technique.

In the present paper, we investigate the differential Krull
dimension of an extension by differentially transcendental elements
in the ordinary case. Even this case has not been
well-understood yet.

Precisely, we study possible bounds for the differential Krull
dimension of a differential polynomial extension. The study was
inspired by works of Jaffard and Seidenberg on the Krull dimension
of polynomial extensions. It should be noted that it is more
difficult to manipulate with differential Krull dimension in
comparison with Krull dimension. In particular, this is because
differential Krull dimension appears together with differential type
and we have to deal with infinite sequences of differential prime
ideals. Therefore, the classic methods used in the dimension theory
of commutative rings do not work in the differential case.

Our main results are Theorem~\ref{HtSCT} and its applications
(Theorems~\ref{dfcheight} and~\ref{S3prop} and their corollaries).
The theorem is a differential analogue of Jaffard's Special Chain
Theorem~\cite[Th\'eor\`eme~3]{Jaf}, which is a very important result
on the behavior of prime chains in an arbitrary polynomial ring (for
applications, see~\cite{AG}, \cite{BHMR}, and~\cite{Jaf}).

Rings with the well-behaved dimension of extensions are called
Jaffard rings, that is rings of finite Krull dimension with $\dim
R[x_1, \ldots, x_n] = R + n$ for all $n$, i.e., the minimal possible
value of the dimension holds. We call the corresponding class of
differential rings, with the minimal possible dimension of
extensions, as J-rings. As an application of Theorem~\ref{HtSCT}, we
obtained that all Jaffard rings are J-rings if they contain rational
numbers or, more generally, are standard (Definition~\ref{stdef})
(this condition is in some sense necessary, as Example~\ref{examp_p}
shows). Also, we proved that, under the same conditions, a
differential ring with a locally nilpotent derivation is also a
J-ring. This shows that both the commutative structure of a
differential ring and properties of the derivation have an impact on
the chains of differential prime ideals of differential polynomial
extensions.

The paper is organized as follows. In Section~\ref{sec1}, we give
definitions and prove basic properties of differential Krull
dimension. Section~\ref{sec2} is devoted to our main result, that is
the Special Chain Theorem (Theorem~\ref{HtSCT}). All technical
lemmas are proved in Section~\ref{sec21}. In Section~\ref{sec22}, we
prove Theorem~\ref{HtSCT} and its immediate corollaries. In
Section~\ref{sec3}, we introduce some classes of differential rings
with ``good'' properties of differential Krull dimension.
Section~\ref{sec31} deals with $\Delta$-arithmetical rings; in
Section~\ref{sec32}, we study locally nilpotent derivations, and, in
Section~\ref{sec33}, we apply the latter material and the Special
Chain Theorem to well-known classes of rings.

\section{Basic definitions and properties}\label{sec1}

All rings considered in this paper are assumed to be commutative and
to contain an identity element. A differential ring is a ring with
finitely many pairwise commuting derivations. An ordinary
differential ring is a ring with one distinguished derivation. An
ideal of a ring is called differential if, for each element of the
ideal, its derivation belongs to the ideal. The set of all prime
differential ideals of a ring $R$ will be denoted by $\diffspec R$
and will be called a differential spectrum.

A homomorphism $f\colon A\to B$, where $A$ and $B$ are differential
rings, is called a differential homomorphism if it commutes with
derivations. For an arbitrary differential ideal $\mathfrak
b\subseteq B$, the ideal $f^{-1}(\mathfrak b)$ is differential and
is called the contraction of $\mathfrak b$. The contraction of
$\mathfrak b$ is denoted by $\mathfrak b^c$. For an arbitrary
differential ideal $\mathfrak a\subseteq A$, the ideal $f(\mathfrak
a)B$ generated by the image of $\mathfrak a$ is differential and is
called the extension of $\mathfrak a$. The extension of $\mathfrak
a$ is denoted by $\mathfrak a^e$. If $\mathfrak q\subseteq B$ is a
differential prime ideal, then its contraction is a prime
differential ideal. Therefore, we have the map $f^*\colon \diffspec
B\to \diffspec A$ by $\mathfrak q\mapsto f^{-1}(\mathfrak
q)=\mathfrak q^c$. For a prime ideal $\mathfrak p\subseteq A$, the
set $(f^*)^{-1}(\mathfrak p)$ is called the fiber over $\mathfrak
p$.

The definition of differential Krull dimension is given by Joseph
Johnson in~\cite{J1}. We will use a slightly modified version of the
definition. Namely, we want the differential type of a differential
ring with the only differential prime ideal to be zero.

\begin{definition}
Let $R$ be a differential ring and $X \subset \diffspec R$. There is
a unique way to define the function called the ``gap''
measure
$$\mu\colon \{\,(\mathfrak p, \mathfrak p')\mid\mathfrak p, \mathfrak p' \in X, \mathfrak p \supset \mathfrak p'\,\}
\to Z \cup {\infty},$$
such that the following conditions hold:
\begin{itemize}
\item[1)] $\mu (\mathfrak p, \mathfrak p') \geq 0$
\item[2)] $\mu (\mathfrak p, \mathfrak p') = 0$ if and only if either $\mathfrak p = \mathfrak
p'$, or there is no infinite descending chain of distinct
differential primes in $X$ such that $\mathfrak p = \mathfrak p_0
\supset \mathfrak p_1 \supset \ldots \supset \mathfrak p'$.
\item[3)] If $d > 0$ is a natural number, $\mu (\mathfrak p, \mathfrak p') \geq d$ if and only if
$\mathfrak p \neq \mathfrak p'$ and there exists an infinite
descending sequence $(\mathfrak p_i)_{i =0, 1, \ldots}$
of differential primes in $X$ such that $\mathfrak p = \mathfrak p_0
\supset \mathfrak p_1 \supset \ldots \supset \mathfrak p'$ and $\mu
(\mathfrak p_{i-1}, \mathfrak p_i) \geq d-1$ for $i = 1, 2, \ldots$.
\end{itemize}

For an arbitrary subset $X\subseteq \diffspec R$, define $\type X$
to be the least upper bound of all the $\mu (\mathfrak p, \mathfrak
p')$, where $\mathfrak p\supseteq \mathfrak p'$ are elements of $X$,
and call it the differential type of $X$. The differential type of a
differential ring $R$ is the differential type of its differential
spectrum and is denoted by $\type R$. We define the differential
dimension of $X$ to be the least upper bound of $n \in Z$ such that
there exists a descending chain $\mathfrak p_0 \supset \mathfrak p_1
\supset \ldots \supset \mathfrak p_n$ of differential prime ideals
in $X$ with $\mu (\mathfrak p_{i-1}, \mathfrak p_i) = \type R$ for
$i = 1, \ldots, n$ and will denote it by $\diffdim X$.
\end{definition}

It should be noted that, in the definition of the differential
dimension of a subset $X$, we use $\type R$ instead of $\type X$ for
pairs of differential prime ideals. If $\diffspec R$ is non-empty,
$\type R$ is well-defined.

\begin{definition}
Let $\mathfrak p, \mathfrak q \in \diffspec R$, $\mathfrak p \subset
\mathfrak q$. Then we will use the following notation:
\begin{itemize}

\item The differential dimension of $\diffspec R$ is
called the differential dimension of $R$ and is denoted by $\diffdim
R$.

\item The differential height of $\mathfrak q$ is
the differential dimension of the set
$$
\{\mathfrak q' \in \diffspec R: \mathfrak q' \subset \mathfrak q\}
$$
and is denoted by $\diffht \mathfrak q$.

\item The differential height of $\mathfrak q$ over $\mathfrak p$
is the differential dimension of the set
$$
\{\mathfrak q' \in \diffspec R: \mathfrak p \subset \mathfrak q'
\subset \mathfrak q\}
$$
and is denoted by  $\diffht \mathfrak q/\mathfrak p$.

\item The differential type of $\mathfrak q$ over $\mathfrak p$ is
the differential type of the set
$$
\{\mathfrak q' \in \diffspec R: \mathfrak p \subset \mathfrak q'
\subset \mathfrak q\}
$$
and is denoted by $\type \mathfrak q/\mathfrak p$.

\end{itemize}
\end{definition}

\begin{remark}\label{remtypedim}
It should be noted that, if $0 < \type R < \infty$, then $\diffdim R
> 0$.
\end{remark}

If $R$ is a differential ring and $n$ is a positive integer, then
$\poly {R}{n}$ will denote the differential polynomial ring $R\{z_1,
\ldots, z_n\}$ in $n$ differential indeterminates over $R$ and, if
$\mathfrak a$ is an ideal of $R$, then $\poly {\mathfrak a}{n}$ will
denote the extension $\mathfrak a\poly{R}{n}$, that is the set of
differential polynomials with coefficients in $\mathfrak a$.

\begin{definition}
Let $\mathfrak p$ be a differential prime ideal of $\poly{R}{n}$ and
$\mathfrak q$ be a differential prime in $R$. We call  $\mathfrak p$
a differential upper to $\mathfrak q$ if an only if $\mathfrak p
\cap R = \mathfrak q$ and $\mathfrak p \neq \poly {\mathfrak q}{n}$.
Otherwise, if $\mathfrak p = \poly {\mathfrak q}{n}$, then
$\mathfrak p$ is called an extended differential prime.
\end{definition}

\begin{definition}
For every prime $\mathfrak p$ in $R$, the field $R_{\mathfrak
p}/\mathfrak pR_{\mathfrak p}$ will be called the residue field of
$\mathfrak p$ and will be denoted by $k(\mathfrak p)$.
\end{definition}

The following theorem proved by Johnson in~\cite[Theorem
of~\S~2]{J1} is the core result on differential Krull dimension that
bears the rest of this section.

\begin{theorem}[Johnson]\label{johnson}
Let $\mathbb K$ be a differential field of characteristic zero with
$m$ derivations. Then, for every $n \in \mathbb N$,
\begin{gather*}
\type \poly{\mathbb K}{n} = m\\
\diffdim \poly{\mathbb K}{n} = n
\end{gather*}
\end{theorem}

To study a fiber over a differential prime ideal via this theorem, we have to require
that the residue field of the differential prime ideal has characteristic zero.

\begin{definition}\label{stdef}
A differential ring $R$ is called standard if, for every
differential prime ideal $\mathfrak p$ of $R$, $k(\mathfrak p)$ has
characteristic zero.
\end{definition}

If $R$ is standard, then, for every $n \ge 1$, $R\{z_1, \ldots,
z_n\}$ is also standard.

\begin{corollary}\label{fiber}
Let $R$ be a standard ring with $m$ derivations. Then, for every
differential prime $\mathfrak p$ of $R$, the fiber over $\mathfrak
p$ for the contraction map has differential type $m$ and dimension
$n$.
\end{corollary}

Ritt algebras are common and the most important examples of
standard rings. However, there are standard rings that are not Ritt
algebras.

\begin{example}\label{example}
Let $R$ be the following ordinary differential ring
$$
\mathbb Z_{(p)}\{x,y,z,u\}/[x^p-py,z^p-x_1-pu,z_1-1].
$$
We claim that $R$ is not the zero ring and, for all prime
differential ideals of $R$, their contractions to $\mathbb Z_{(p)}$
are zero.

To prove the first assertion, we consider $\mathbb Q\otimes_{\mathbb
Z_{(p)}} R$. Hence, we have
\begin{gather*}
\mathbb Q\otimes_{\mathbb Z_{(p)}} R = \mathbb Q\{x,z,y,u\}/[z_1-1,
y-x^p/p, u- (z^2-x_1)/p]=\\
=\mathbb Q\{x,z\}/[z_1-1]=\mathbb Q\{x\}[z].
\end{gather*}
In particular, we have just computed the fiber over the zero ideal,
that coincides with the differential spectrum of the ring $R$.

To prove the second assertion, we consider $\mathbb F_p
\otimes_{\mathbb Z_{(p)}} R$. Hence, we have
\begin{multline*}
\mathbb F_p\otimes_{\mathbb Z_{(p)}} R = \mathbb
F_p\{x,y,u,z\}/[x^p,z^p-x_1,z_1-1]=\\
=\mathbb F_p
\{y,u\}[x,x_1,z,z_1]/(x^p, z^p-x_1, z_1-1)=\mathbb
F_p\{y,u\}[z,x]/(x^p).
\end{multline*}
Let $J$ denote the smallest radical differential ideal of this
ring. We will show that $J$ is the whole ring. Indeed, since
$x^p=0$, $x\in J$. Since $J$ is differential, $x_1\in J$. From the
equality $z^p=x_1$, it follows that $z\in J$. But $z_1=1$. So, $1\in
J$.
\end{example}

It should be noted that, if we are interested in the differential
spectrum of a standard differential ring $R$, then we can suppose
that $R$ is a Ritt algebra. Indeed, we can replace $R$ by the ring
$\mathbb Q\otimes_{\mathbb Z}R$. The differential spectrum of the
latter ring coincides with the differential spectrum of the ring
$R$.

Since differential rings of characteristic $p$ usually have big
spectra, the assertion of Theorem~\ref{johnson} might be wrong, as the
following example shows.

\begin{example}\label{examp_p}
Let $L=\mathbb F_p(x_n)_{n\in \mathbb N}$ be a purely transcendental
extension of the finite field $\mathbb F_p$, the derivation
$\partial$ is set to be zero. Then the ring of differential
polynomials $L\{y\}$ has infinite differential type.

Indeed, for every subset $X\subseteq \mathbb N$, we define $I_X$ to
be the ideal generated by $\{y_k^p-x_k\mid k\in X\}$. The ideals
$I_X$ are differential prime ideals. In $\mathbb N$, one can find a
descending sequence of subsets $X_k$ such that $X_k\setminus
X_{k+1}$ is countable. Thus, the family $I_{X_k}$ is a descending
chain of differential prime ideals. Since each set $X_k\setminus
X_{k+1}$ is infinite, we can find a descending chain of sets $Y_{k,
n}$ such that $Y_{k,n}\setminus Y_{k, n+1}$ is countable. Then we
have a descending chain of the corresponding differential prime
ideals $I_{Y_{k,n}}$.

Repeating this procedure for each $Y_{k,n}$, we can produce a chain
of differential prime ideals of any differential type. Thus, the
differential type of $L\{y\}$ is infinite.
\end{example}

\begin{remark}\label{comdim}
Let $S$ be a ring of finite Krull dimension, then we have the inequality~(\cite[Theorem~2]{S1})
$$\dim S + 1 \leq \dim S[x] \leq 2\cdot \dim S + 1.$$
\end{remark}

The next lemma shows possible bounds for the differential type of
$R\{z\}$ in the spirit of the inequality in Remark~\ref{comdim}.

\begin{lemma} \label{typeineq}
Let $R$ be a standard differential ring of differential type $t$
with $m$ derivations. Then we have
$$
\max (t, m) \le \type R\{z_1,\ldots, z_n\} \le t + m.
$$
\end{lemma}
\begin{proof}
Since every chain in $R$ extends to a chain in $R\{z_1, \ldots,
z_n\}$, then $t \le \type \poly{R}{n}$. Furthermore, there is a
chain of differential type $m$ by Corollary~\ref{fiber}.

Let $\mathfrak C$ be a chain realizing the differential type of
$\poly{R}{n}$ and $\mathfrak C \cap R$ be the chain in $R$
consisting of contractions of ideals of $\mathfrak C$. Since there
is no chain of differential type $m + 1$ in the fiber over any
differential prime, we have $\type (\mathfrak C \cap R) \geq \type
\poly{R}{n} - m$. Therefore, $\type R \geq \type \poly{R}{n} - m$.
\end{proof}

\begin{corollary}\label{growtype}
Let $R$ be a standard differential rings of differential type $t$
with $m$ derivations and $n$ is a positive integer. Then
$$
\type \poly{R}{n} \leq \type \poly{R}{n+1}.
$$
\end{corollary}
\begin{proof}
Since $m \leq \type \poly{R}{n}$ by Lemma~\ref{typeineq}, then the
inequality is clear after regarding $\poly{R}{n+1} =
\poly{R}{n}\{z_{n+1}\}$.
\end{proof}

\begin{corollary}\label{chaintype}
Let $R$ be a standard differential ring with $m$ derivations. 
If $\mathfrak p, \mathfrak q \in \diffspec R$, $\mathfrak p \subset \mathfrak q$ 
and $\type \mathfrak q/\mathfrak p = t$. Then, for any $\mathfrak q'$ upper to
$\mathfrak q$ in $\poly{R}{n}$, we have $\type \mathfrak
q'/\poly{\mathfrak p}{n} \le t + m$.
\end{corollary}
\begin{proof}
Since the ring $R_{\mathfrak q}/{\mathfrak p}R_{\mathfrak q}$ has differential type $t$,
the result follows from Lemma~\ref{typeineq}.
\end{proof}

\begin{lemma}\label{gendim}
Let $R$ be a standard differential ring of differential type $t$
such that $\type \poly {R}{n} = \max (t, m)$. Then
\begin{itemize}
    \item[a)] $n \leq \diffdim \poly {R}{n}$ if $t < m$;
    \item[b)] $n + \diffdim R \leq \diffdim \poly {R}{n}$ if $t = m$;
    \item[c)] $\diffdim R \leq \diffdim \poly{R}{n}$ if $t > m$.
\end{itemize}
\end{lemma}
\begin{proof}
In the case of $\diffdim R = \infty$ the last two assertions are trivial,
since every chain in $R$ can be extended to $\poly{R}{n}$.
The case of $t < m$ is straight-forwardly deduced from Corollary~\ref{fiber}.

Furthermore, let $\mathfrak p_0 \subset \ldots, \subset
\mathfrak p_d$ be a chain realizing the differential dimension of $R$.
Then, by Corollary~\ref{fiber}, there is a chain $\poly{\mathfrak p_0}{n}
\subset \ldots \subset \poly{\mathfrak p_d}{n} \subset \mathfrak q_1
\subset \ldots \subset \mathfrak q_n$, where all $\mathfrak q_i$ are
uppers to $\mathfrak p_d$. If $t = m$, this chain has dimension at least $d + n = \diffdim R + n$.
In the last case, the inequality holds since every chain in $R$ can be extended to $\poly {R}{n}$
(the type of $\mathfrak q_1 \subset \ldots \subset \mathfrak q_n$ less then $\type \poly{R}{n} = t$,
so the chain has dimension zero).
\end{proof}

There is no reason for dimension to be bounded above in general, but
in special cases upper bound exists.

\begin{lemma}\label{specdim}
Let $R$ be a standard differential ring of differential type $t$
with $m$ derivations. Then
\begin{enumerate}
    \item $n \leq \diffdim \poly {R}{n} \leq n\cdot\diffdim R$ if $t =
    0$;
    \item $1 \leq \diffdim \poly{R}{n} \leq \diffdim R$ if $\type \poly {R}{n} = t + m$.
\end{enumerate}
\end{lemma}
\begin{proof}

(1). By Lemma~\ref{typeineq}, the type of $\poly{R}{n}$ is $m$. By
Corollary~\ref{fiber}, the fiber over every differential prime ideal
has dimension $n$, therefore, $\diffdim \poly{R}{n} \leq n \cdot
\diffdim R$. And the left-hand part follows Lemma~\ref{gendim}.

(2). The left-hand inequality is obvious. Since the contraction of a
chain of differential type $t + m$ in $\poly {R}{n}$ is a
differential type $t$ chain by Corollary~\ref{chaintype}, the
right-hand inequality is clear.
\end{proof}

\begin{definition}\label{Jdef}
Let $R$ be a differential ring of finite differential type with $m$ derivations. 
We say that $R$ is a J-ring if, for every $n > 0$, $\max (t, m) = \type
R\{z_1, \ldots, z_n\}$ and
\begin{enumerate}
    \item $\diffdim R\{z_1, \ldots, z_n\} = n$ if $t < m$;
    \item $\diffdim R < \infty$ and $\diffdim R\{z_1, \ldots, z_n\} = n + \diffdim R$ if $t = m$;
    \item $\diffdim R < \infty$ and $\diffdim R\{z_1, \ldots, z_n\} = \diffdim R$ if $t > m$.
\end{enumerate}
\end{definition}

These rings have the minimal possible values for both the
differential type and the differential dimension of $\poly{R}{n}$,
i.e., the lower bound of the inequalities in Lemma~\ref{gendim}
holds. This is a differential analogue of Jaffard rings, that is,
rings $R$ with $\dim R[x_1, \ldots, x_n] = \dim R + n$ for every
$n$, i.e., the lower bound of the inequality in Remark~\ref{comdim}
holds. The class of Jaffard rings is relatively wide, for example,
finite dimensional Noetherian rings are Jaffard. Our goal is to
establish similar results for the class of J-rings.

\begin{lemma}\label{infhtprop}
Let $R$ be a differential ring of finite differential type $t$
and $\mathfrak p$ be a differential prime ideal in $R$.
If $\diffht \mathfrak p = \infty$, then, for every integer $n \geq 0$,
there exists a differential prime $\mathfrak p_n \subset \mathfrak p$
such that  $\diffht \mathfrak p_n = n$.
\end{lemma}
\begin{proof}
Suppose that the contrary holds. Obviously, if
there is a differential prime ideal $\mathfrak q$ with $\diffht \mathfrak
q = n$, then, for every $k\leqslant n$, there is a differential prime
ideal $\mathfrak q$ with $\diffht \mathfrak q = k$. Now, we may suppose that
there exists $k$ such that our claim is false for $k$, then, for
every differential prime $\mathfrak q \subset \mathfrak p$ such that
$\diffht \mathfrak q \geq k$, $\diffht \mathfrak q = \infty$.

Let $\mathfrak q$ be a differential prime of infinite height in $R$.
We will find a differential prime $\mathfrak q'$ in $R$ such that
$\type \mathfrak q/\mathfrak q' = t$ and $\diffht \mathfrak q' = \infty$.
Since $\diffht \mathfrak q = \infty$, there is a chain
$$\mathfrak {p}^{k+1} \subset \ldots \subset \mathfrak {p}^0 \subset \mathfrak q$$
such that $\type\mathfrak {p}^i/\mathfrak {p}^{i+1} = t$.
Since $\diffht \mathfrak p^1 \geq k$,  $\diffht \mathfrak p^1 = \infty$.
Also $t = \type \mathfrak p^0/\mathfrak p^1 \leq \type \mathfrak q/\mathfrak p^1 \leq t$,
therefore $\mathfrak q' = \mathfrak p^1$ is the required ideal.

Let $\mathfrak q_0 = \mathfrak p$. Applying the procedure above to $\mathfrak q_i$ and
denoting $\mathfrak q_{i+1} = \mathfrak q_i'$, we will obtain the chain
$\mathfrak p = \mathfrak q_0 \supset \mathfrak q_1 \supset \mathfrak q_2 \supset \ldots$ such
that $\type \mathfrak q_i/\mathfrak q_{i+1} = t$. Hence, we have $\type R > t$, a contradiction.
\end{proof}

The next two lemmas are quite similar, the only difference arises from the definition of J-rings. 
Namely, the case of differential rings with differential type less then the number 
of derivations is slightly different from another ones; 
the finiteness of the differential dimension is not assumed.

\begin{lemma}\label{facloclt}
Let $R$ be a differential ring of finite differential type $t < m$.
If $R$ is not a J-ring, then, for at least one pair of differential
prime ideals $\mathfrak p \subset \mathfrak m$, $R_{\mathfrak
m}/{\mathfrak p}R_{\mathfrak m}$ is not a J-ring.
\end{lemma}
\begin{proof}
Suppose that the definition of a J-ring does not hold for a natural
number $n > 0$. If $\type \poly{R}{n} > \max(t,m)$, then let
$\mathfrak q_1 \subset \mathfrak q_2$ be differential primes in
$\poly{R}{n}$ such that $\type \mathfrak q_2/\mathfrak q_1 = \type
\poly{R}{n}$. Then $\mathfrak q_1 \cap R \subset \mathfrak q_2 \cap
R$ is the required pair.

Now, suppose that $\type \poly{R}{n}=m=\max(t,m)$. Since
$t<m$, we only need to show that the first condition of
Definition~\ref{Jdef} does not hold. Consider a chain $\mathfrak q_1 \subset
\ldots \subset \mathfrak q_r$ in $\poly{R}{n}$ such that $r =
\diffdim \poly{R}{n}$ if $\diffdim \poly{R}{n} < \infty$, or $r = n
+ 1$ otherwise. In both cases, we have $r>n$. Then,
obviously, our assertion holds for $\mathfrak q_1 \cap R \subset
\mathfrak q_r \cap R$.
\end{proof}

\begin{lemma}\label{facloc}
Let $R$ be a differential ring of finite differential type $t \geq
m$. Suppose that the ring $R$ satisfies the following property: if
$\type \poly{R}{n} = t$ for every $n$, then $\diffdim R < \infty$.
If $R$ is not a J-ring, then, for at least one pair of differential
prime ideals $\mathfrak p \subset \mathfrak m$, $R_{\mathfrak
m}/{\mathfrak p}R_{\mathfrak m}$ is not a J-ring.
\end{lemma}
\begin{proof}
Suppose that, for some natural number $n > 0$, $\type \poly{R}{n} >
t = \max(t,m)$. Let $\mathfrak q_1 \subset \mathfrak q_2$ be
differential primes in $\poly{R}{n}$ such that $\type \mathfrak
q_2/\mathfrak q_1 = \type \poly{R}{n}$. Then $\mathfrak q_1 \cap R
\subset \mathfrak q_2 \cap R$  is the required pair.

Now, suppose that $\type \poly{R}{n} = t$ for every $n$.
Since $t\geq m$, we should cotradict to the second and the third
conditions of Definition~\ref{Jdef}. Under the hypothesis of the
proposition, $\diffdim R < \infty$. Let $n$ be a natural number such
that the definition of J-ring does not hold for $n$. There are two
cases, either $\diffdim \poly{R}{n}$ is finite or not. If
$\diffdim\poly{R}{n} = \infty$, there is a chain $\mathfrak q_1
\subset \ldots \subset \mathfrak q_r$ with $r=\diffdim R + n + 1$. Otherwise, we
have either $t = m$ and $\diffdim \poly{R}{n} > n+\diffdim R$, or $t > m$
and $\diffdim\poly{R}{n} > \diffdim R$. In these cases, we consider a
chain $\mathfrak q_1 \subset \ldots \subset \mathfrak q_r$ in
$\poly{R}{n}$ such that $r = \diffdim \poly{R}{n}$. In all cases, we
see that $\mathfrak q_1 \cap R \subset \mathfrak q_r \cap R$ is the
desired pair of prime differential ideals of $R$.
\end{proof}

Since we are mostly dealing with ordinary differential rings in the
rest of the paper, the following lemma will be very useful.

\begin{lemma}\label{niceht}
Let $\mathbb K$ be an ordinary differential field of characteristic
zero. Then, for every nonzero differential prime ideal $\mathfrak p$
of $\mathbb K\{z\}$,
$$
\diffht \mathfrak p = 1.
$$

Moreover, for every infinite descending chain $C$ of prime
differential ideals, we have
$$
\bigcap_{\mathfrak q \in C} \mathfrak q = 0.
$$
\end{lemma}
\begin{proof}

By Proposition~3 of~\cite[Chapter~III, Section~2]{Ros}, we obtain
the first assertion. Since the algebra $K\{z\}/\mathfrak p$ is
differentially algebraic over $K$, the differential type of this
quotient is $0$ by Theorem in Section~2 of~\cite{J1}. Hence, all
descending chains of differential prime ideals of $K\{z\}/\mathfrak p$
are finite.

\end{proof}

\begin{remark}\label {whydown}
By Ritt-Raudenbush's basis theorem~\cite[Chapter II, Section 4, Theorem 1]{K1}, 
all ascending chains in $\poly{\mathbb K}{n}$ have finite length if $\mathbb K$ has characteristic zero, 
therefore, an infinite chain has to be descending.
\end{remark}

\section{Special Chain Theorem}\label{sec2}
In this section we prove our main result, that is
Theorem~\ref{HtSCT}. As an immediate corollary, we obtain that all
ordinary standard Jaffard rings are J-rings.

\subsection{Auxiliary lemmas}\label{sec21}

\begin{lemma}\label{type0}
Let $R$ be an ordinary standard differential ring with the condition
$\type R\{z\} = 1$. Let $\mathfrak q' \subset \mathfrak q$ be
differential primes in $R$ such that $\type \mathfrak q/\mathfrak q'
= 0$, $\mathfrak p'\subseteq \mathfrak p$ be differential primes in
$R\{z\}$ of finite differential height such that $\mathfrak p' \cap
R = \mathfrak q'$ and $\mathfrak p$ is a differential upper to
$\mathfrak q$. If $\diffht \mathfrak p = \diffht \mathfrak p' + 1$,
then $\diffht \mathfrak p = \diffht \mathfrak q\{z\} + 1$.
\end{lemma}

\begin{proof}
The inequality $\diffht \mathfrak p \geq \diffht \mathfrak q\{z\} +
1$ follows from Theorem~\ref{johnson}.

Since there is a chain $\mathfrak C$ between $\mathfrak p'$ and
$\mathfrak p$ of differential type $1$ and all descending chains
between $\mathfrak q'$ and $\mathfrak q$ have a finite length, there
is a differential prime $\mathfrak t$ in $R$ between $\mathfrak q'$
and $\mathfrak q$ such that there are infinite number of elements of
$\mathfrak C$ lying over $\mathfrak t$. By Remark~\ref{whydown}, the
chain over $\mathfrak t$ has differential type $1$, and,
consequently, differential dimension $1$, so $\diffht \mathfrak
p/\mathfrak t\{z\} = 1$. Therefore, $\diffht \mathfrak t\{z\} \leq
\diffht\mathfrak p - 1$.

Since $\mathfrak p' \subset \mathfrak p_\alpha$ for each $\mathfrak
p_\alpha \in \mathfrak C$ and
$$
\mathfrak t\{z\} = \bigcap_{\substack{\mathfrak p_\alpha \in \mathfrak C\\
\mathfrak p_\alpha \cap R = \mathfrak t}} \mathfrak p_\alpha
$$
by Lemma~\ref{niceht}, we have $\mathfrak p' \subset \mathfrak
t\{z\}$ and, therefore, $\diffht \mathfrak p' \leq \diffht \mathfrak
t\{z\}$.

So, we have
$$
\diffht\mathfrak p = \diffht\mathfrak p'+ 1\leq \diffht \mathfrak t\{z\}
+ 1 \leq \diffht \mathfrak q\{z\} + 1.
$$

\end{proof}

\begin{lemma}\label{htlemma}
Let $R$ be an ordinary standard differential ring of finite differential type,
$n \geq 1$ be an integer such that $\type \poly{R}{n} = \type R$,
and $\mathfrak p$ be a differential prime ideal in $\poly{R}{n}$.
If $\diffht \mathfrak p < \infty$, then $\diffht (\mathfrak p \cap R) < \infty$.
\end{lemma}
\begin{proof}
Since every chain of differential prime ideals in $R$ can be
extended to a chain in $\poly{R}{n}$ of the same differential type,
the differential height of the contraction of $\mathfrak p$ cannot
be infinite.
\end{proof}

\begin{corollary}\label{fincon}
Let $R$ be a standard differential ring with one derivation such
that $\type R \leq 1$. Let $n \ge 1$ be an integer such that $\type
\poly {R}{n} = 1$. If a differential prime $\mathfrak p$ in
$\poly{R}{n}$ has finite differential height, then, for every $0 < k
< n$, $\mathfrak p \cap \poly{R}{k}$ also has finite differential
height.
\end{corollary}
\begin{proof}
By Corollary~\ref{growtype}, we have $\type \poly{R}{k} = 1$
for $k>0$. Therefore, the result follows from
Lemma~\ref{htlemma}.
\end{proof}

\begin{lemma}\label{ml}

Let $R$ be an ordinary standard differential ring such that
$$
\type R \le 1 \mbox{ and } \type R\{z\} = 1.
$$
Let $\mathfrak q$ be a differential prime ideal in $R$ and
$\mathfrak p$ be a differential upper to $\mathfrak q$ in $R\{z\}$.
If $\diffht p < \infty$, then
$$
\diffht \mathfrak p = \diffht \mathfrak q\{z\} + 1.
$$
If we additionally suppose that, for some natural number $n > 0$,
$\type \poly{R}{n} = 1$ and $\mathfrak p\{z_2,\ldots,z_n\}$ has
finite differential height, then
$$
\diffht \mathfrak p\{z_2, \ldots, z_n\} = \diffht \mathfrak q\{z_1,
\ldots, z_n\} + 1.
$$
\end{lemma}
\begin{proof}
The inequality $\diffht \mathfrak p \ge \diffht \mathfrak q\{z\} +
1$ follows from Lemma~\ref{niceht}. So, we have to show the inverse
one.

Since $\diffht \mathfrak p/\mathfrak q\{z\} = 1$, then $m=\diffht
\mathfrak p \ge 1$. Let $\mathfrak p_1$ be a differential prime with
$\diffht \mathfrak p_1 = m - 1$ in a chain realizing differential
height of $\mathfrak p$ and $\mathfrak q_1 = \mathfrak p_1 \cap R$.
If $\type R = 0$, then $\type \mathfrak q/\mathfrak q_1 = 0$ and we
are done by Lemma~\ref{type0}. So, we may suppose that $\type R =
1$.

We use the induction by $\diffht \mathfrak q$. If $\diffht \mathfrak
q = 0$, then $\type \mathfrak q/\mathfrak q_1=0$ and we are done by
Lemma~\ref{type0}.

If $\diffht \mathfrak q/\mathfrak q_1 \geq 1$, then
$$
\diffht \mathfrak p = \diffht \mathfrak p_1 + 1 \leq \diffht \mathfrak q_1\{z\} + 2 \leq \diffht \mathfrak q\{z\} + 1.
$$

As for the second assertion, we have
$$
\mathfrak p\{z_2, \ldots, z_n\} \cap R\{z_2, \ldots, z_n\} =
\mathfrak q\{z_2, \ldots, z_n\}
$$
and, after considering the ring $R\{z_2, \ldots, z_n\}$ as $R$,
$\mathfrak p\{z_2, \ldots, z_n\}$ as $\mathfrak p$, and $\mathfrak
q\{z_2,\ldots,z_n\}$ as $\mathfrak q$, the result follows from the
first part of the lemma.

\end{proof}

\begin{proposition}\label{smalltype}
Let $R$ be an ordinary standard differential ring such that $\type R
\leq 1$, $n \ge 1$ be an integer number such that $\type \poly{R}{n}
= 1$, and $\mathfrak p$ be a differential upper to $\mathfrak q$ in
$\poly {R}{n}$. If $\mathfrak p$ has a finite differential height,
then
$$\diffht \mathfrak p = \diffht \poly{\mathfrak q}{n} + \diffht
\mathfrak p/\poly{\mathfrak q}{n}.$$
\end{proposition}
\begin{proof}
The inequality $\diffht \mathfrak p \ge \diffht \poly{\mathfrak q}{n}
+ \diffht \mathfrak p/\poly{\mathfrak q}{n}$
is clear. So, we must show the inverse one.

We induce by $n$. The case $n = 1$ immediately
follows from Lemma~\ref{ml}. Therefore, we assume that the result
holds for all $k < n$. Set $\mathfrak p_1 = \mathfrak p \cap
R\{z_1\}$, then $\mathfrak p_1$ has a finite differential height by
Corollary~\ref{fincon}. If $\mathfrak p_1 = \mathfrak q\{z_1\}$,
then, representing $R\{z_1, \ldots, z_n\}$ as
$R\{z_1\}\{z_2,\ldots,z_n\}$, we derive the result from the
induction hypothesis.

If $\mathfrak q\{z_1\} \subset \mathfrak p_1$, then Lemma \ref{ml}
implies that
$$
\diffht \mathfrak p_1\{z_2, \ldots, z_n\} = \diffht \mathfrak
q\{z_1, \ldots, z_n\} + 1
$$
and if $\mathfrak p = \mathfrak p_1\{z_2,\ldots,z_n\}$, then we are
done. Now, we may suppose that $\mathfrak p \supset \mathfrak p_1\{z_2,\ldots,z_n\}$.
Therefore, the induction assumption implies that
\begin{align*}
\diffht \mathfrak p &= \diffht \mathfrak p_1\{z_2, \ldots, z_n\} +
\diffht \mathfrak p/\mathfrak p_1\{z_2, \ldots, z_n\} = \\
&=1 + \diffht \mathfrak q\{z_1, \ldots, z_n\} + \diffht \mathfrak
p/\mathfrak p_1\{z_2, \ldots, z_n\}.
\end{align*}
Since $\mathfrak q\{z_1\} \subset \mathfrak p_1$,  $\diffht
\mathfrak p_1/\mathfrak q\{z_1\} = 1$ by Lemma~\ref{niceht}. Also,
since every differential prime $\mathfrak q'$ in $R\{z_1\}$ can be
extended to the differential prime $\mathfrak q'\{z_2,\ldots,z_n\}$
in $\poly{R}{n}$, every chain between $\mathfrak q\{z_1\}$ and
$\mathfrak p_1$ can be extended to a chain between $\mathfrak
q\{z_1, \ldots, z_n\}$ and  $\mathfrak p_1\{z_2, \dots, z_n\}$.
Therefore, we have
$$
1 + \diffht \mathfrak p/\mathfrak p_1\{z_2, \ldots, z_n\} \le
\diffht \mathfrak p/\mathfrak q\{z_1, \ldots, z_n\}.
$$
The latter inequality implies that
$$
\diffht \mathfrak p
\le \diffht \mathfrak q\{z_1, \ldots, z_n\} + \diffht \mathfrak
p/\mathfrak q\{z_1, \ldots, z_n\}.
$$
\end{proof}

\begin{lemma}\label{typejump}
Let $R$ be an ordinary standard differential ring satisfying the
condition $\type \poly{R}{n} = t \geq 2$. Let $\mathfrak q' \subset
\mathfrak q$ be differential primes in $R$ such that $\type
\mathfrak q/\mathfrak q' < t$, $\mathfrak p' \subseteq \mathfrak p$
be differential primes in $\poly{R}{n}$ of finite differential
height such that $\mathfrak p' \cap R = \mathfrak q'$ and $\mathfrak
p$ is a differential upper to $\mathfrak q$. If $\diffht \mathfrak p
= \diffht \mathfrak p' + 1$, then $\diffht \mathfrak p = \diffht
\poly{\mathfrak q}{n}$.
\end{lemma}
\begin{proof}
The inequality $\diffht \mathfrak p \geq \diffht \poly{\mathfrak
q}{n}$ is clear, we must show the inverse one.

By Corollary~\ref{chaintype} and the assumption on $\type \mathfrak
q'/\mathfrak q$, we have $\type \mathfrak q/\mathfrak q' = t - 1$.
There is the only way to get a chain of differential type $t$
between $\mathfrak p$ and $\mathfrak p'$. Namely, there is a chain
$C = \{\mathfrak q' \subset \ldots \subset \mathfrak q_{\alpha}
\subset \ldots \subset \mathfrak q_0 \subset \mathfrak q \}$ such
that, for every $\mathfrak q_{\alpha} \in C$, there exists
$\mathfrak p_{\alpha}$ being a differential upper to $\mathfrak
q_{\alpha}$ such that $\mathfrak p_{\alpha} \subset \poly{\mathfrak
q_{\alpha - 1}}{n}$ for every $\alpha$. Therefore, $\diffht
\mathfrak p = \diffht \poly{\mathfrak q_0}{n} \leq \diffht \poly
{\mathfrak q}{n}$.
\end{proof}

\begin{lemma}\label{bigtype1in}
Let $R$ be an ordinary standard differential ring with the condition
$\type R = t> 1.$ Let $\mathfrak q$ be a differential prime in $R$
and $\mathfrak p$ be an upper to $\mathfrak q$ in $R\{z\}$ of finite
differential height. Then
$$
\diffht \mathfrak p = \diffht \mathfrak q\{z\},
$$
and, for any positive integer $n$ such that $\diffht \mathfrak
p\{z_2,\ldots,z_n\}<\infty$, we have
$$
\diffht \mathfrak p\{z_2, \ldots, z_n\} = \diffht \poly{\mathfrak q}{n}.
$$
\end{lemma}
\begin{proof}
The inequality $\diffht \mathfrak p \ge \diffht \mathfrak q\{z\}$ is
clear. If $\diffht \mathfrak p = 0$, we are done. So, we may suppose
that $\diffht \mathfrak p \ge 1$. Let $\mathfrak p'$ be a
differential prime in a chain realizing the differential height of
$\mathfrak p$ with $\diffht \mathfrak p' = \diffht \mathfrak p - 1$
and let $\mathfrak q' = \mathfrak p' \cap R$. From
Corollary~\ref{chaintype}, it follows that

$$
\type R\{z\} - 1 \le \type \mathfrak q/\mathfrak q' \le \type R\le
\type R\{z\}.
$$

The case of $\type \mathfrak q/\mathfrak q' = \type R\{z\} - 1$
follows from Lemma~\ref{typejump}.

In the second case, by Lemma~\ref{htlemma}, we have that $\diffht
\mathfrak q < \infty$ and $\diffht \mathfrak q' < \diffht \mathfrak
q$.  We will use the induction by $\diffht \mathfrak q$. In the case
of $\diffht \mathfrak q = 1$, we have $\diffht \mathfrak p' =
\diffht \mathfrak q'\{z\}$ by the first part of the proof.
Otherwise, by the induction hypothesis, we have $\diffht \mathfrak
p' = \diffht \mathfrak q'\{z\}$. Therefore,

$$
\diffht \mathfrak p = \diffht \mathfrak p' + 1 =
\diffht \mathfrak q'\{z\} + 1 \le \diffht \mathfrak q\{z\}.
$$

Since $\mathfrak q\{z_2, \ldots, z_n\} = R\{z_2, \ldots, z_n\} \cap
\mathfrak p\{z_2, \ldots, z_n\}$, then, after considering
$\poly{R}{n}$ as $R\{z_2, \ldots, z_n\}\{z_1\}$, the second
assertion follows from the first one.
\end{proof}

\begin{lemma}\label{bigtype}
Let $R$ be an ordinary standard differential ring satisfying the
condition $\type R = t \ge 1$. Let $\mathfrak p$ be a differential
prime of finite differential height in $\poly{R}{n}$ and $\mathfrak
p \cap R = \mathfrak q$. Then, for every positive integer $n$ such
that $\type \poly{R}{n} \ge 2$, we have $\diffht \mathfrak p =
\diffht \poly{\mathfrak q}{n}$.
\end{lemma}
\begin{proof}
If $\mathfrak p = \poly{\mathfrak q}{n}$ our claim is trivial.

The case of $\type R = 1$ follows from Lemma~\ref{typejump} since
$\type \poly{R}{n} \ge 2$. We may suppose that $\type R > 1$, then,
by Lemma~\ref{typeineq}, $\type \poly{R}{n} \geq \type R > 1$ for
every $n$.

We use the induction by $n$. The base follows from
Lemma~\ref{bigtype1in}. Therefore, we assume that the result holds
for all $k < n$. Set $\mathfrak p_1 = \mathfrak p \cap
\poly{R}{n-1}$. By Lemma~\ref{bigtype1in}, $\diffht \mathfrak
p_1\{z_n\} = \diffht \mathfrak p$, so $\mathfrak p_1\{z_n\}$ has
finite differential height. If $\mathfrak p_1\{z_n\}=\poly{\mathfrak
q}{n}$, then we are done. Otherwise, $\mathfrak p_1\{z_n\}$ is an
upper to $\mathfrak q\{z_n\}$ because $\mathfrak p_1\{z_n\}\cap
R\{z_n\} = \mathfrak q\{z_n\}$. After regarding $\poly{R}{n}$ as
$\poly{R\{z_n\}}{n-1}$, we have $\diffht \mathfrak p_1\{z_n\} =
\diffht \poly{\mathfrak q}{n}$ by the induction hypothesis.
\end{proof}

\subsection{The main theorem and its corollaries}\label{sec22}

\begin{theorem}[Special Chain Theorem]\label{HtSCT}
Let $R$ be an ordinary standard differential ring of finite differential type,
then, for every positive integer $n$, the following holds
$$
\diffht \mathfrak p = \diffht \poly{\mathfrak q}{n} + \diffht
\mathfrak p/\poly{\mathfrak q}{n},
$$
where $\mathfrak p \in \diffspec \poly{R}{n}$ and $\mathfrak p \cap R = \mathfrak q$.
\end{theorem}
\begin{proof}
Let $\diffht \mathfrak p < \infty$. If $\type \poly{R}{n} = 1$, the
result follows from Proposition~\ref{smalltype}, otherwise from
Lemma~\ref{bigtype}.

If $\diffht \mathfrak p = \infty$, then it follows from
Lemma~\ref{infhtprop} that, for every $i$, there is $\mathfrak p_i
\subset \mathfrak p$ with $\diffht \mathfrak p_i = i$.
By the case of finite differential
height, we have $\diffht \poly{(\mathfrak p_i\cap R)}{n} \geq i - n$.
Therefore, $\diffht \poly{\mathfrak q}{n} = \infty$.
\end{proof}

This differential analogue of Jaffard's theorem (see~\cite[2.3]{AG}
and~\cite[Theorem~1]{BHMR}) will be our main tool for the
investigation of the differential dimension of certain classes of
rings.

Jaffard proved in~\cite[Th\'eor\`eme~2]{Jaf} that the sequence
$$
d_n = \dim R[x_1, \ldots, x_n] - \dim R[x_1, \ldots, x_{n-1}]
$$
is eventually constant and the eventual value $d_R$ is at most $\dim
R + 1$. Therefore, there is an equality $\dim R[x_1, \ldots, x_n] =
d_R\cdot n + C$ for some constant $C$ and sufficiently large $n$-th.

For each ideal $\mathfrak q \subset \poly{R}{n}$, we use the
following notation
$$
\mathfrak q_m = \mathfrak q \cap R[z_1, \ldots,
\partial^{m} z_1, \ldots, z_n, \ldots,
\partial^{m}z_n].
$$

Let $R$ be an ordinary standard differential ring of finite Krull dimension,
then $\type R = 0$ and, by Lemma~\ref{typeineq}, $\type \poly{R}{n}
= 1$ for every $n$. Also, by Lemma~\ref{specdim}, $\diffdim
\poly{R}{n} < \infty$.

\begin{theorem}\label{dfcheight}
Let $R$ be an ordinary standard differential ring of finite Krull
dimension. If, for some $n$, we have $\diffdim \poly{R}{n} > n$,
then $d_R > 1$.
\end{theorem}
\begin{proof}
Let $\mathfrak m$ be a differential prime in $\poly{R}{n}$ and
$\mathfrak q = \mathfrak m \cap R$ such that $\diffht \mathfrak m =
\diffdim \poly {R}{n} > n$. By Theorem~\ref{HtSCT}, we have $\diffht
\mathfrak m = \diffht \poly {\mathfrak q}{n} + \diffht \mathfrak
m/\poly{\mathfrak q}{n}$. Also, by Corollary~\ref{fiber}, $\diffht
\mathfrak m/\poly{\mathfrak q}{n} \leq n$. Therefore, $\diffht
\poly{\mathfrak q}{n} > 0$ and, by the argument similar to the proof
of Lemma~\ref{type0}, there are a differential prime $\mathfrak t
\subset \mathfrak q$ in $R$ and a differential prime ideal
$\mathfrak p$ in $\poly{R}{n}$ that is upper to $\mathfrak
t$ such that $\mathfrak p \subset \poly {\mathfrak
q}{n}$. Since we have the inclusion
\begin{align*}
\mathfrak p_m\subset\
 \mathfrak q[z_1, \ldots, \partial^{m} z_1, \ldots, z_n,
\ldots, \partial^{m}z_n]
\end{align*}
and there is a chain of $n\cdot(m + 1)$ primes in
$R[z_1, \ldots, \partial^{m} z_1, \ldots, z_n, \ldots, \partial^{m}z_n]$ over $\mathfrak q$,
then $\height{\mathfrak p_m} + n\cdot (m + 1) < d_R\cdot n\cdot (m+1) + C$.
Assume that $d_R = 1$, then we have $\height \mathfrak p_m < C$.

There is an infinite chain $\mathfrak p = \mathfrak p^0 \supset
\mathfrak p^1 \supset \ldots \supset \poly {\mathfrak t}{n}$ of
differential type $1$. Let $f_i$ be polynomials such that $f_i \in
\mathfrak p^{i-1}$ and $f_i\notin \mathfrak p^i$. Suppose $k_j$ is
the minimal integer such that
$$
f_i \in R[z_1, \ldots, \partial^{k_j} z_1, \ldots, z_n, \ldots, \partial^{k_j}z_n] \mbox{ for all } i < j.
$$
Therefore, there is a chain $\mathfrak p_{k_j} \supset \mathfrak
p^1_{k_j} \supset \ldots \supset \mathfrak p^j_{k_j}$ and
$$
\height \mathfrak p_{k_j} > j.
$$
A contradiction with $d_R = 1$.
\end{proof}

\begin{corollary}\label{jaffard}
If $R$ is a Jaffard ring that is ordinary and standard, then $R$ is
a J-ring.
\end{corollary}
\begin{proof}
By definition, $R$ has finite Krull dimension and $d_n = 1$ for
every $n$.
\end{proof}

\begin{corollary}
If $R$ is a Noetherian standard differential ring with one
derivation, then $R$ is a J-ring.
\end{corollary}
\begin{proof}
Since a Noetherian local ring has finite Krull dimension
(\cite[Corollary~11.11]{AM}), all chains in $R$ have finite length.
Hence, $\type R = 0$.

Supoe that our assumption is wrong. Since $\type R = 0$ and $R$ is ordinary,
by Lemma~\ref{facloclt},
$R_{\mathfrak m}$ is not a J-ring for some $\mathfrak m$. But
$R_{\mathfrak m}$ is local Noetherian and, hence, Jaffard
by~\cite[Corollary~2]{BHMR}.
\end{proof}

\begin{remark}
In the case of $\type R < \type \poly{R}{n}$, one should not expect
a ``good'' behavior of differential height. For example, let $R$ be
the Nagata example~\cite[A1, p.~203, Example~1]{Nag} of a Noetherian
ring with $\dim R = \infty$ and $\partial$ be the zero derivation.
We have $\type R = 0$, and there are prime ideals such that $\diffht
\mathfrak q = k$ for every $k$. However, for every prime ideal
$\mathfrak q$ in $R$, $\diffht \mathfrak q\{z\} = 0$. Furthermore,
we expect that there exists a ring with $\diffht \mathfrak p =
\infty$ and $\diffht \mathfrak p\{z\} = 0$.
\end{remark}

\section{Examples of J-rings}\label{sec3}
In this section we are going to establish some well-known classes of
rings included in the class of J-rings.

\subsection{$\Delta$-arithmetical rings}\label{sec31}
Our first class is an analogue of arithmetical rings. Arithmetical
rings are important since every Pr\"ufer domain is
arithmetical.
\begin{definition}
A ring $A$ is called arithmetical if every finitely generated ideal
is locally principal, that is, for every finitely generated ideal
$\mathfrak a$, the image of $\mathfrak a$ in $A_{\mathfrak p}$ is
principal for every prime ideal $\mathfrak p$ in $A$.
\end{definition}
Arithmetical rings are Jaffard~\cite[Theorem~4]{S2}.

\begin{definition}
A $\Delta$-arithmetical ring is a differential ring $A$
such that, for every finitely generated ideal $\mathfrak a$, $\mathfrak a
A_{\mathfrak p}$ is principle ideal for every differential prime
$\mathfrak p$.
\end{definition}

\begin{remark}
It is well-known that an integer domain is arithmetical if and only
if it is a Pr\"ufer domain\cite[Chapter~1,
Section~6, Theorem~62]{Kap}. Consider a local $\Delta$-arithmetical
domain, such that the maximal ideal is differential. Then this ring
is a local Pr\"ufer domain, i.e. a valuation domain.
\end{remark}

The following example shows that there are $\Delta$-arithmetical
rings that are not arithmetical.

\begin{example}
Let $O$ be the Krull example of an integrally closed
one-dimen\-sional local domain not being a valuation
ring~\cite[p.~670f]{Kr}. We will recall the construction. Let $K$ be
an algebraically closed field, $x$ and $y$ are indeterminates. The
ring $O$ consists of the rational functions $r(x,y)$ such that $x$
does not divide denominator of $r(x,y)$ and $r(0,y) \in K$. Let $y$
and elements of $K$ be constant and $x' = 1$.

The ring $O$ is local and is not a valuation ring, so, $O$ is not a
Pr\"ufer domain. The only differential prime is $(0)$ and, in the
field of fractions of $O$, every ideal is principle.
\end{example}

\begin{proposition}\label{pruferuz}
Let $\mathfrak p$ be a differential prime ideal of a
$\Delta$-arithmetical domain $O$. If $\mathfrak q$ is a differential
upper to zero, then $\mathfrak q \nsubseteq \mathfrak p\{z_1, \ldots, z_n\}$.
\end{proposition}
\begin{proof}
We will follow Seidenberg's proof~(\cite[Theorem~4]{S1}).

We may assume that $O$ is local and $\mathfrak p$ is the maximal
ideal. Then $O$ is a valuation domain. Let $f \in \mathfrak q$, then
$f = c\cdot g$, where $c \in \mathfrak p$ and $g$
has at least one coefficient equal to 1.
Therefore, $g \in \mathfrak q$ and $g \notin \mathfrak p\{z_1,
\ldots, z_n\}$.
\end{proof}

\subsection{Derivations with special properties}\label{sec32}

Let us recall that the characteristic of a ring $R$ is the natural
number $n$ such that $n\mathbb Z$ is the kernel of the unique ring
homomorphism from $Z$ to $R$.

\begin{lemma}\label{uzconst}
Let $O$ be a differential integral domain of characteristic zero,
$\mathfrak p$ is a differential upper to zero in
$O\{z_1, \ldots, z_n\}$. If there is a nonzero polynomial $f \in
\mathfrak p$ such that its coefficients are constant for some
$\partial \in \Delta$, then there is $f^* \in \mathfrak p$ such that
$f^*$ has relatively prime coefficients in $\mathbb Z$.
\end{lemma}
\begin{proof}
Every non-zero polynomial $f$ with constant coefficients
can be presented as the following sum
$$
a_1f_1 + \ldots + a_nf_n,
$$
where $a_i\in O$ are constant and $f_i$ are polynomials with integer
coefficients. Let $f\in \mathfrak p$ be taken such that $n$ is minimal
among all such numbers. If $n=1$, then $f=af_1\in \mathfrak p$. Since
$\mathfrak p$ is prime and $\mathfrak p\cap O=0$, $f_1\in \mathfrak p$.

Suppose that $n>1$, then
$$
\partial f = a_1\partial f_1 + \ldots + a_n\partial f_n.
$$
Set $g$ to be the following polynomial
$$
f_1\partial f - \partial f_1f = a_2g_2 + \ldots + a_ng_n,
$$
where $g_i = \partial f_if_1 - \partial f_1f_i$. The polynomial $g$
has constant coefficients. From the definition of $f$, it follows
that $g=0$. Let $L$ be the field of fractions of $O$, then the
fraction $f_1/f$ is a $\partial$-constant in $L\langle
y_1,\ldots,y_n\rangle$. Therefore, $f_1/f$ belongs to the subfield
of $\partial$-constants of $L$~\cite[Chapter~II, Section~9,
Corollary~5 of Theorem~4]{K1}. In particular, $af_1 = b f$ for some
elements $a,b\in O$. Thus, $af_1\in \mathfrak p$ and again $f_1\in
\mathfrak p$.

So, we have a polynomial $f$ with integral coefficients in $\mathfrak
p$. We may suppose that coefficients of $f$ are coprime because
$f=df'$, where $d$ is the greatest common divisor of the
coefficients. Since $\mathfrak p$ is an upper to zero, $d\notin \mathfrak
p$. Thus, $f'\in \mathfrak p$.
\end{proof}

\begin{remark}\label{ckalg}
If the ring $O$ is a Ritt algebra in the previous theorem, then all
coefficients of $f^*$ are invertible.
\end{remark}

\begin{lemma}\label{lnd}
Let $(O, \mathfrak m)$ be a differential local integral domain of
characteristic zero and $\mathfrak p$ be an upper to zero
in $\diffspec O\{z_1, \ldots, z_n\}$. If there is $f \in \mathfrak p$
with the property: there is $\partial\in \Delta$ such that every coefficient $a$ of $f$
satisfies $\partial^n(a)=  0$ for some natural $n$.
Then there exists $h \in \mathfrak p$ having relatively
prime coefficients in $\mathbb Z$.
\end{lemma}
\begin{proof}
We will say that $a\in O$ is $\partial$-nilpotent if
$\partial^n(a)=0$ for some $n$. If $a\in O$ is $\partial$-nilpotent
than we define $d(a)$ to be the maximal $n\in \mathbb N$ such that
$\partial^n(a)\neq 0$.

Every polynomial $f$ with $\partial$-nilpotent coefficients can be
presented as the following sum
$$
a_1f_1+\ldots+a_nf_n,
$$
where $a_i\in O$ are $\partial$-nilpotent, $d(a_1)\geqslant
d(a_2)\geqslant\ldots\geqslant d(a_n)$, and all $f_i$ have integer
coefficients. For such representation of $f$, we set $d(f)=d(a_1)$
and define $n(f)$ to be the maximal number $i$ such that
$d(a_i)=d(f)$.

Let $f$ be a non-zero polynomial of $\mathfrak p$ with
$\partial$-nilpotent coefficients and its representation
$$
a_1f_1+\ldots+a_nf_n
$$
be chosen such that the pair $(d(f),n(f))$ is lexicographically
minimal. If $d(f)=0$, then the lemma follows from
Lemma~\ref{uzconst}. Suppose that $d(f) > 0$. Let $g$ be as follows
$$
g=f_1\partial f - \partial f_1f = a_2g_2 + \ldots + a_ng_n +
\partial(a_1)h_1+\ldots+\partial(a_n)h_n,
$$
where $g_i = \partial f_if_1 - \partial f_1f_i$, $h_i=f_1f_i$. Then
$g$ has $\partial$-nilpotent coefficients and this representation
for $g$ is less than the initial representation of $f$. Thus, $g=0$.
So, we have $f_1\partial f = \partial f_1 f$. As in the proof of
Lemma~\ref{uzconst}, we derive that $af_1\in \mathfrak p$ for some
$a\in O$. And since $\mathfrak p\cap O = 0$, $f_1\in \mathfrak p$.

\end{proof}

\begin{lemma}\label{lnduz}
Let $(O, \mathfrak m)$ be a differential local integral domain of
characteristic zero, $\mathfrak p$ is a differential upper to zero
in $O\{z_1, \ldots, z_n\}$. If there is $f \in \mathfrak p$ such
that all coefficients are $\partial$-nilpotent for some $\partial\in
\Delta$, then $\mathfrak p \nsubseteq \mathfrak m\{z_1, \ldots, z_n\}$.
\end{lemma}
\begin{proof}
By Lemma~\ref{lnd}, there is $f^* \in \mathfrak p$ with coefficients
in $\mathbb Z$. If $\mathfrak m \cap \mathbb Z = 0$ then $O \supset
\mathbb Q$ and, by Remark \ref{ckalg}, $f^* \notin \poly{\mathfrak
m}{n}$. Otherwise, since $O$ is local, $\mathbb Z_{(p)} \subset O$
for some prime $p$. Therefore, since coefficients are relatively
prime, at least one coefficient is not divisible by $p$. Thus, $f^*
\notin \poly{\mathfrak m}{n}$.
\end{proof}

\subsection{An application of the Special Chain
Theorem}\label{sec33}

We will say that a differential ring $R$ satisfies the property {\bf
S3}, if, for every natural number $n > 0$, for every differential
prime ideals $\mathfrak q_1 \subset \mathfrak q_2$ in $R$, there is
no differential prime $\mathfrak p$ in $\poly{R}{n}$ such that
$\mathfrak p$ is an upper to $\mathfrak q_1$ and $\mathfrak p
\subset \poly{\mathfrak q_2}{n}$. This property has a sense like
being a Stable Strong S-Domain for usual polynomial extensions.
\begin{remark}\label{s3}
If {\bf S3} holds, every element of a chain that ends in
$\poly{\mathfrak q}{n}$ is an extended prime too. Therefore, for
every chain $\mathfrak C$ that ends in extended prime, $\type
\mathfrak C = \type \mathfrak C \cap R$.
\end{remark}

\begin{theorem}\label{S3prop}
Let $R$ be an ordinary standard differential ring of finite
differential type that satisfies {\bf S3}. Suppose that either
$\type R = 0$ or $\diffdim R < \infty$, then $R$ is a J-ring.
\end{theorem}
\begin{proof}
Let $n > 0$ be a natural number. By Lemma~\ref{typeineq}, $1 \leq
\type \poly {R}{n} < \infty$. First, we will show that $\type \poly
{R}{n} = \max(\type R, 1)$.

If there exists such differential prime $\mathfrak q$ in $R$ such
that $\diffht \poly{\mathfrak q}{n} > 0$, then applying
Remark~\ref{s3} to the chain realizing differential height of
$\mathfrak q$, we have $\type R \geq \type \poly{R}{n}$. Therefore,
$\type R = \type \poly{R}{n}$.

In the second case, all extended primes have differential height
zero. Suppose that $\type R = \type \poly{R}{n}$ then,
by Remark~\ref{remtypedim}, there exists $\mathfrak q \subset R$
of nonzero differential height, so $\diffht \poly{\mathfrak q}{n} > 0$.
Hence $\type R < \type \poly{R}{n}$. By Remark~\ref{remtypedim}, there exist
$\mathfrak p \subset \poly{R}{n}$ of nonzero differential height. By
Theorem~\ref{HtSCT}, we have
$$
0 < \diffht \mathfrak p = \diffht \poly{(\mathfrak p\cap R)}{n}
+ \diffht \mathfrak p/\poly{(\mathfrak p \cap R)}{n} = \diffht
\mathfrak p/\poly{(\mathfrak p \cap R)}{n}.
$$

But by Theorem~\ref{johnson}, $\type \mathfrak p/\poly{\mathfrak p \cap
R}{n} = 1$. Therefore, $\type \poly{R}{n} = 1 = \max(\type R, 1)$.

Now, we will prove the assertion about the differential dimension.
If $\type R = 0$, then $\type \poly{R}{n} = 1$ and the property {\bf
S3} implies that, for every extended prime, $\diffht \poly{\mathfrak
q}{n} = 0$. By Theorem~\ref{HtSCT}, $\diffht \mathfrak p = \diffht
\mathfrak p/\poly{\mathfrak p \cap R}{n} \leq n$ for every
differential prime $\mathfrak p$ and, by Theorem~\ref{johnson},
there exists $\mathfrak p$ such that $\diffht \mathfrak
p/\poly{\mathfrak p \cap R}{n} = n$. So, $\diffdim \poly {R}{n} =
n$.

By Theorem~\ref{HtSCT}, $\diffht \mathfrak p =
\diffht \poly{(\mathfrak p\cap R)}{n} + \diffht p/\poly{(\mathfrak p\cap R)}{n}$
for every differential prime $\mathfrak p$ in $\poly{R}{n}$.
Since {\bf S3} holds, $\diffht \poly{(\mathfrak p\cap R)}{n} = \diffht \mathfrak p\cap R$.
In the case of $\type R = 1$, we have $\diffht \mathfrak p \leq \diffdim R + n$.
Consequently, $\diffdim \poly {R}{n} = \diffdim R + n$.
Otherwise, if $\type R > 1$, we have $\diffht \mathfrak p = \diffht \poly{(\mathfrak p\cap R)}{n} \leq \diffdim R$.
Hence, $\diffdim \poly {R}{n} = \diffdim R$.
\end{proof}

\begin{remark}
Since valuation domains of finite Krull dimension are
Jaffard(\cite[Corollary~2]{BHMR}, \cite[Theorem~4]{S2}), we could
use Corollary~\ref{jaffard} to prove that all $\Delta$-arithmetical
rings of finite Krull dimension are J-rings. However, we are able to
prove this without the assumption to have finite Krull dimension.
\end{remark}

\begin{corollary}
Let $R$ be an ordinary standard differential ring of finite
differential type that has either differential type zero or finite
differential dimension. If $R$ is $\Delta$-arithmetical then $R$ is
a J-ring.
\end{corollary}
\begin{proof}
Let $n > 0$ be an integer, $\mathfrak q' \subset \mathfrak q$
differential primes in $R$ and $\mathfrak p$ be a differential upper
to $\mathfrak q'$ in $\poly{R}{n}$. An application of
Proposition~\ref{pruferuz} to $R/\mathfrak q'$ shows that $R$
satisfies {\bf S3}. By the previous theorem, $R$ is J-ring.
\end{proof}

\begin{corollary}
Let $R$ be an ordinary standard differential ring of finite differential type
with locally nilpotent derivation that has either differential type zero or
finite differential dimension. Then $R$ is a J-ring.
\end{corollary}
\begin{proof}
Let $n > 0$ be an integer, $\mathfrak q' \subset \mathfrak q$
differential primes in $R$ and $\mathfrak p$ be a differential upper
to $\mathfrak q'$ in $\poly{R}{n}$. Suppose that $\mathfrak p
\subset \poly{\mathfrak q}{n}$. Let $f \in \mathfrak pR_{\mathfrak
q}/\mathfrak q'R_{\mathfrak q}$, then $c\cdot f$ has
$\partial$-nilpotent coefficients, where $c$ is the product of coefficients denominators of $f$.
Since $R_{\mathfrak q}/\mathfrak q'R_{\mathfrak q}$
is a local integral domain, we have a contradiction with
Lemma~\ref{lnduz}. Therefore, $R$ satisfies {\bf S3} and, consequently, is a J-ring.
\end{proof}

We have shown that chains of differential prime ideals in polynomial
extensions are affected by both the commutative structure of the
base ring (Jaffard rings and $\Delta$-arithmetical rings) and by
properties of derivations (locally nilpotent derivation). This
implies that standard rings that are not J-rings have to be rare.
Namely, an example of a standard ordinary differential ring
of finite Krull dimension that is not a J-ring should be a non-Jaffard ring that
admits non-trivial derivation. Moreover, an example of a standard
ordinary ring such that $\max (\type R, 1) < \type R\{z\}$
have to be even more complicated.

\section{Acknowledgement}
The author is very grateful to Evgeny Golod, who taught him
commutative algebra and have acquainted him with the works of Gilmer
and Seidenberg. The author is deeply appreciative to Alexander Levin for his advise and support.
The author is very grateful to Dmitry Trushin, who taught him differential algebra and
gave many crucial and useful suggestions including examples~\ref{example} and~\ref{examp_p}.

\bibliographystyle{model1-num-names}
\bibliography{diffkdim}

\end{document}